\documentclass[a4paper,10pt]{article}

\usepackage[colorlinks, linkcolor=blue,  citecolor=blue, urlcolor=blue]{hyperref}
\usepackage[colorlinks]{hyperref}
\usepackage{amsfonts}
\usepackage{amssymb}
\usepackage{amsmath}
\usepackage{amsthm}
\usepackage{enumerate}
\usepackage{epsfig}
\usepackage{mathrsfs}
\usepackage{graphicx}
\usepackage{dsfont}
\usepackage{epsfig,color}
\usepackage{dsfont}
\usepackage[utf8]{inputenc}
\usepackage{bm}
\newtheorem{definicion}{Definition}[section]
\newtheorem{definition}[definicion]{Definition}
\newtheorem{proposition}[definicion]{Proposition}
\newtheorem{theorem}[definicion]{Theorem}

\newtheorem{lemma}[definicion]{Lemma}
\newtheorem{remark}[definicion]{Remark}

\numberwithin{equation}{section}

\def\fin { \vskip 0pt \hfill \hbox{\vrule height 5pt width 5pt depth 0pt} \vskip 12pt}

\newcommand{\Pw}[2]{\mathscr{P}_{\overrightarrow{#1}^{#2}}(\Omega)}
\newcommand{\dmi}[3]{W^2_{#1}(#2,#3)}
\newcommand{\dm}[3]{\mathbf{W}^2_{\overrightarrow{#1}}(\overrightarrow{#2},\overrightarrow{#3})}
\newcommand{\dmk}[4]{\mathbf{W}^2_{\overrightarrow{#1}^{#4+1}}(\overrightarrow{#2}^{#4},\overrightarrow{#3}^{#4+1})}
\begin{document}

\title{{\large\textbf{On the Convergence of the JKO-scheme and Blow-up of solutions for a Multi-species Chemotaxis System with no Mass Preservation}}}
\date{}

\author{\textsc{\small Julio C. Valencia-Guevara}\\
	\text{\small Departamento de Matemática y Estadística,}\\
	\text{\small Universidad Católica San Pablo, Arequipa, Perú}\\
	{\small julioguevara08@gmail.com, jcvalencia@ucsp.edu.pe}}

\maketitle

%\footnotetext[1]{Email addresses: julioguevara08@gmail.com, jcvalencia@ucsp.edu.pe}
%\footnotetext{\textbf{MSC 2010.} {\it Primary}: 92C17; {\it Secondary}: 35K40, 35K55, 35K45, 35B44, 58Exx}
%\footnotetext{\textbf{Keywords}: Blow-up; Chemotaxis; Keller-Segel; Multispecies; JKO scheme; Optimal Transport.}

\begin{abstract}
This work considers a chemotaxis system for multi-species that includes birth or death rate terms, which implies no mass preservation of the populations. We aim to show the convergence to a $L^{\infty}\ -\ $weak solutions, that is local in time, of the JKO - scheme arising from the Optimal Transport Theory, in the spirit of \cite{petrelli-tudorasco,carrillo-santambrogio}. Currently, $L^{\infty}$ solutions have shown to be important in order to get uniqueness. Since death rate case does not ensure global solutions, for arbitrary initial data, in this framework, it could be interest to analyze the Blowing-up phenomenon of this system. Therefore, in the last section, we get sufficient conditions that implies blowing-up phenomenon in finite time and we draw several stages where this occurs. This last part can be seen as a partial generalization of the blowing-up results in \cite{conca-espejo-vilches}.  \newline
	
	\noindent\textbf{MSC 2010.} {\it Primary}: 92C17; {\it Secondary}: 35K40, 35K55, 35K45, 35B44, 58Exx,
	28A33 \vspace{0.1cm}
	
	\noindent\textbf{Keywords}: Blow-up; Chemotaxis; Keller-Segel; Multispecies; JKO scheme; Optimal Transport.
	
\end{abstract}

%============================================================================================================
%=======================================================INTRODUCTION=========================================

\section{Introduction}

Mathematical models in chemotaxis have been and continue to be studied widely. Chemotaxis is a mechanism by which some organisms move in response to a chemical signal produced by a substance. Historically, the pioneering works in mathematical models about chemotaxis were conducted by Patlak \cite{patlak} and Keller $\&$ Segel \cite{keller-segel}. In a simplified version that model can be written as:

\begin{equation}
\left\{
\begin{array}{rcll}
\displaystyle\frac{\partial \rho}{\partial t}&=&\displaystyle\Delta\rho-\alpha\nabla \cdot\left(\rho\nabla v\right),&\text{ for }(x,t)\in \Omega\times (0,T),\\
-\Delta v&=&\rho,&\text{ for }(x,t)\in \Omega\times (0,T),\\
\rho(x,0)&=&\rho_{0}(x)\geq0&\text{ for }x\in\Omega,
\end{array}
\right.\label{eq:keller-segel}
\end{equation}
under suitable boundary conditions, where $\rho$ denotes the population density of the organism, $v$ the concentration of the chemical substance, $\alpha$ the sensitivity of the bacteria to the chemo-attractant and  $\Omega\subset\mathbb{R}^2$.

As it was pointed out in the monography \cite{blanchet} the model proposed by Patlak-Keller-Segel has been an important mathematical ingredient in the study of biological processes such as: pattern formation of cells for meiosis, embryo-genesis, etc. The reader interested in more details about these applications can also consult the references in \cite{blanchet}. Despite the fact that \eqref{eq:keller-segel} models the dynamic of a single specie reacting to a single chemical signal, there is rich mathematical literature in which system \eqref{eq:keller-segel} and other related are analysed. A good revision of mathematical models in chemotaxis can be found in \cite{hillen-painter} and references therein.

In \cite{jager-luckhaus}, it was shown the existence of solutions, global in time, under a smallness condition over the parameter $\alpha>0$ and the initial mass $m_0=\int_{\Omega}\rho_0(x)\ dx$. In this same work, it was shown that radially symmetric solutions blow-up when the term $\alpha m_0$ is greater than a critical value that was not specified. Nowadays, such critical value, for $m_0=1$, is well known to be $8\pi$. More precisely, in \cite{dolbeault-perthame,blanchet-dolbeault} the following dichotomy was shown: {\it a)} if $\alpha m_0<8\pi$ (sub-critical case) then there is a global in time weak solution (this is a classic solution for regular bounded domains), {\it b)} if $\alpha m_0>8\pi$ (super-critical case) then, blowing-up of the solutions occurs in a finite time. A study for the critical case, that is when $\alpha m_0=8\pi$, was conducted by Blanchet et al. in \cite{blanchet-carrillo-masmoudi}. In their work the authors showed the existence of solutions global in time and which concentrate mass for $t\to\infty$ to a delta Dirac. Also for the critical case, in \cite{blanchet-carlen-carrillo}  some qualitative properties such as asymptotic behaviour of solutions and existence of basins of attraction were studied.
 
On the other hand, Optimal Transport theory has been used to analyse system \eqref{eq:keller-segel} as it is evidenced in the works \cite{blanchet-calvez-carrillo,carrillo-santambrogio}. In \cite{blanchet-calvez-carrillo} it was shown the convergence of the JKO-scheme to a global in time weak solution and in \cite{carrillo-santambrogio} it was performed a $L^{\infty}$ bound for the JKO-scheme and its convergence to a local in time weak solution of \eqref{eq:keller-segel} was shown. The results in \cite{carrillo-santambrogio} were obtained for a non linear diffusion. Optimal Transport theory has also allowed to obtain some results about uniqueness  as in \cite{carrillo-lisini-mainini}, where the authors reached this result by showing that integrable and $L^{\infty}$ solutions satisfy a gradient flow formulation of a suitable convex functional. Also in \cite{liu-jinhuan} a degenerate Keller-Segel model is treated by building bounded solutions and showing uniqueness of $L^{\infty}$ weak solutions. These works indicate that bounded solutions are essential in order to get uniqueness results at least for the moment. In other hand, numerical methods based in the JKO scheme (or gradient flow scheme) are being recently developed, in particular for the simplest version of Keller-Segel model the reader can consult \cite{carrillo-mathes-wolfram}. 

In this work we are concerned with the multispecies system for chemotaxis:
 \begin{equation}
 \left\{
 \begin{array}{rcll}
 \displaystyle\frac{\partial \rho_{i}}{\partial t}&=&\displaystyle\Delta\rho_{i}-\alpha_i\nabla \cdot\left(\rho_{i}\nabla v\right)+H_i(\overrightarrow{\rho}),&\text{ in }\text{Int}(\Omega)\times (0,T),\\
 -\Delta v&=&\sum_{j=1}^{N}\alpha_j\rho_j,&\text{ in }\text{Int}(\Omega)\times (0,T),\\
 (\alpha_i\rho_i\nabla v-\nabla\rho_i)\cdot \overrightarrow{n}&=&0,\ v=0&\text{ in } \partial\Omega\times (0,T),\\
 \rho_i(x,0)&=&\rho_{0,i}(x)\geq0&\text{ for }x\in\Omega,
 \end{array}
 \right.\label{eq:system-origin}
 \end{equation}
 for $i=1,2,\ldots,N$. Here, $\overrightarrow{\rho}(x,t)=(\rho_1,\ldots,\rho_N)$ denotes the vector of densities $\rho_i\in L^1_+(\Omega)$ of the bacteria population $i$, $\Omega\subset\mathbb{R}^2$ is a sufficient smooth convex bounded set, the parameter $\alpha_i$ represents the sensitivity of the population $i$ to the chemical substance $v$ and $\overrightarrow{H}(\overrightarrow{\rho})=(H_1,H_2,\ldots,H_N)$ denotes the vector representing the birth or death rate of the specie $i$. In this paper, we consider only the case when either all the $H_i$ is birth rate or death rate i.e. $H_i\geq0$ for all $i$ or $H_i\leq0$ for all $i$ respectively.
 
Several models generalizing those started by Keller, Segel and Patlak for multispecies have been proposed. One of the most general can be found in \cite{horstman}. These systems consider for example several populations with nonlinear diffusion coefficients under the effect of various chemical signals also with nonlinear diffusions. Even in \cite{horstman} it was distinguished the case when the chemotaxis models admit a Lyapunov functional, in this work, it is also analysed the linear stability, variational structure, etc. Let us comment that the system \eqref{eq:system-origin} can be seen as a particular case of those generalizations found in \cite{horstman} where the presence of terms $H_i(\overrightarrow{\rho})$ is justified because, from a biological point of view, ``\textit{several of the observed effects might be caused by growth terms and not by the chemical forces''}.
 
In \cite{Kavallaris-ricciardi-zecca} a similar system to \eqref{eq:system-origin} was considered but in absence of the growth terms $H_i(\overrightarrow{\rho})$. The authors consider the cases of production and destruction of the chemical substances, i.e. $\alpha_i>0$ and $\alpha_i<0$ respectively  and they analyse the variational structure by showing the existence of a Lyapunov functional. Then, they establish an infinite dimensional version of the logarithm Hardy-Littlewood-Sobolev inequality for this particular case. Other important studies were conducted in \cite{wolansky-multispecies,espejo-stevens-velazquez,conca-espejo-vilches}. In \cite{wolansky-multispecies} a mutispecies system is considered without growth terms and it is given conditions on the parameters in order to get existence of solutions. Also equilibria states are characterized as critical points of the free energy functional associated. In \cite{espejo-stevens-velazquez} results concerned to the blowing-up of solutions for two species are obtained, showing, among other facts, that blowing-up of one specie can imply the blowing-up of the second. In \cite{conca-espejo-vilches} the picture is complemented for the case of two species. The authors identify a curve into the plane of the masses that replaces the concept critical mass for a single specie. This curve separates the regions where blowing-up occurs and global in time solutions are being reached.
  
The aim of this work is to employ the classic JKO-scheme, arising from the Optimal Transport theory \cite{villani,santambrogio-book}, that was used for first time in the pioneering work by Otto et al. in \cite{JKO}. We show existence of weak solutions that satisfies a local in time $L^{\infty}$ bound and some regularity in $W^{1,2}(\Omega)$ is gained. This technique was widely developed in \cite{ambrosiogiglisavare} where it is applied to solve many evolution equations that can be interpreted as the continuity equation. It is worth mentioning some works from the large list of references that have use such technique. For parabolic-parabolic Keller-Segel system refer to \cite{BL,blanchet-carrillo-kinderlehrer}, for magnetic fluids \cite{Otto}, for a periodic setting \cite{CarrilloSlepcev,VF}, for fractional operators \cite{erbar} among others. The second goal of this work is to analyse the blowing-up of solutions for the case when death rate is present in the model. This is because $L^{\infty}-$ estimates for the JKO-scheme blow-up in finite time no matter the initial mass. Thus, it could be interesting to get conditions that implies such blowing-up. We focus on the model for two species and we take radially symmetric initial conditions over the unitary disk and we take $H_i(\overrightarrow{\rho})=-c_i\rho_i$.

Note that system \eqref{eq:system-origin} can be rewritten formally as:
 
 $$\frac{\partial \rho_i}{\partial t}=\nabla\cdot\left(\rho_i\nabla\frac{\delta\mathcal{F}}{\delta\rho_i}[\overrightarrow{\rho}]\right)+H_i(\overrightarrow{\rho}),\ i=1,2,\ldots,N$$ 
 where $\mathcal{F}$ is the free energy functional defined in \eqref{eq:energy-functional}. This system can be interpreted as a perturbed gradient flow system over the space of positive finite measures on $\Omega$ and with respect to the well known Riemmannian structure induced by Wasserstein metric. This interpretation is merely formal because the Wasserstein metric is defined over the space of finite, positive measures with constant prescribed mass and the  system \eqref{eq:system-origin} does not have the mass conservation property. Thus, we can not expect that the solutions live in a such space. However, as it was performed in \cite{petrelli-tudorasco} the JKO-scheme for gradient flows can still be used to build weak solutions. Let us sketch this idea, formally the system above can be interpreted as:
 
$$u'=-\nabla F(u)+H(u),$$
in an Euclidean setting. Then, based on an implicit Euler scheme, we use the discretization 

\begin{equation}
	\frac{u_{k+1}-u_{k}-\tau H(u_{k})}{\tau}=-\nabla F(u_{k+1}).\label{eq:implicit-euler-scheme}
\end{equation}
The last equality allows to built $u_{k+1}$, for given $u_{k}$, as
$$u_{k+1}\in \text{Arg}\min_{u}\left\{F(u)+\frac{|u-v_{k}|^2}{2\tau}\right\},$$
where $v_{k}=u_k+\tau H(u_{k}).$

This papers follows the ideas of \cite{carrillo-santambrogio,petrelli-tudorasco} which allow to consider mass variations in the right hand of a conservative system as \eqref{eq:keller-segel}. Simultaneously, we can get uniform $L^{\infty}$ estimates for the discrete scheme. Please note  that as far as we are concern, there is no work dealing with decreasing mass term in multispecies models via Optimal Transport theory. Unlike the usual procedure using the JKO-scheme for conservative models, in this case weak convergence is not enough to pass the limit on the discrete version of \eqref{eq:system-origin}. Consequently, by using the $L^{\infty}$ estimates we obtain a $L^{1}$ strong compactness result for the approximate solution that permits circumvent the previous difficulty. In fact, for the spacial variables $L^{2}$ compactness holds true. As a final important complement, we analyze conditions under which blowing-up in finite time occurs. We concentrate in the case of two species and take $H_i(\overrightarrow{\rho})=-c_i\rho_i$, for $c_i>0$. In fact, we draw several stages where this phenomenon happens and generalize some results of \cite{conca-espejo-vilches}. As conducted in \cite{carrillo-santambrogio}, this work is also valid for non-linear diffusions but for simplicity purposes we prefer to keep this study only for linear diffusion. Based in the recent work \cite{carrillo-mathes-wolfram}, the present work opens the possibility of develop a numerical method using the JKO - scheme for multispecies and with no mass preservation.

During the progress of this work we were advised of the recent work by Karmakar and Wolansky \cite{karmakar-wolansky} which may have some similarities to ours. With respect to this fact, let us highlight some important differences. The main difference in the model resides in the growth term $H_i(\overrightarrow{\rho})$ which is not considered in \cite{karmakar-wolansky}. Besides, we can see \cite{karmakar-wolansky} as a generalization for multi-species of the work of Blanchet et al. in \cite{blanchet-calvez-carrillo} that is intended for a single specie. However, \cite{blanchet-calvez-carrillo} has been complemented by the work of Carrillo and Santambrogio in \cite{carrillo-santambrogio}. The present work can be considered as a generalization of \cite{carrillo-santambrogio} in the conservative case ($H_i=0$). Thus, both works  \cite{karmakar-wolansky} and the present complement each other. In addition to that, as it was commented above $L^{\infty}$ solutions play an important role when one wishes to show uniqueness of solutions in chemotaxis models. Hence, our work is developed in a framework that would allow to show this uniqueness result. From a technical point of view, there is also a difference with the treatment of the non-linearity in the model. In \cite{karmakar-wolansky} the non-linearity is treated by using a flow interchange technique taken from the work \cite{mathes-maccann}, while in this work we use regularity result coming from the Poisson equation taken from \cite{Augusto-ponce}. Finally, the last one difference that we consider important is the analysis that we carry out regarding the blowing-up of solutions. Therefore, our work and \cite{karmakar-wolansky} consider important and different frameworks for multispecies chemotaxis systems.

This paper is organized as follows, section \ref{section:preliminaries} provides the mathematical tools and the basic notations that will be used in this work. Section \ref{section:Linfty-estimate} deals with the construction of the JKO-scheme following \cite{carrillo-santambrogio} for multispecies. Section \ref{section:JKO-scheme} provides the $L^{\infty}$ bound in the cases of birth rate and death rate. Section \ref{sectio:convergence} focuses on the compactness and convergence of the JKO-scheme to get a weak solution. Section \ref{section:blow-up}  is devoted to analyze conditions implying blowing-up of a system with two species.
 
 \section{Mathematical Notation and Preliminaries}
 \label{section:preliminaries}
 In this section we provide the mathematical preliminaries and establish the basic notations that will be used in this work. Through this work $\Omega\subset\mathbb{R}^2$ denotes, by simplicity, a compact convex set with non empty interior $\text{Int}(\Omega)$. Next, we denote by $\mathcal{M}^+(\Omega)$ the set of non negative finite Borel measures in $\Omega$. Let $\rho\in\mathcal{M}^+(\Omega)$ so that it is absolutely continuous with respect to the Lebesgue Measure $dx$ then, we use the abuse of notation to represent with $\rho(x)$ the density of $\rho$ with respect to $dx$.
 
 Since we do not have mass preservation property, it does not make sense to work on the classic probability measures space. Therefore, we use the following space that depends on the total mass.
 \begin{definition}
 	\label{defin:wasserstein-vector-space}
 	Let us $\overrightarrow{m}=(m_1,m_2,\ldots,m_N)\in\mathbb{R}^N$ with $m_i>0$, we define the set
 	\begin{equation}
 		\label{eq:defin-wasserstein-vector-space}
 			\Pw{m}{\ }=\left\{\overrightarrow{\rho}=(\rho_1,\rho_2,\ldots,\rho_N): \rho_i\in\mathcal{M}^+(\Omega),\ \rho_i(\Omega)=m_i \right\}.
 	\end{equation}
 \end{definition}
 
 We can endow $\Pw{m}{\ }$ with a metric arising from Optimal Transport theory. Indeed, given $\mu,\rho\in\mathcal{M}^+(\Omega)$ such that $m:=\mu(\Omega)=\rho(\Omega)$ the 2-Wasserstein distance between $\mu$ and $\rho$ is defined as,
 \begin{equation}
 	\label{eq:wasserstein-metric-classic}
 	\dmi{m}{\mu}{\rho}=\inf_{\gamma\in\Gamma(\mu,\rho)}\int_{\Omega\times\Omega}|x-y|^2\ d\gamma(x,y),
 \end{equation}
 where,
 \begin{equation*}
 	\Gamma(\mu,\rho)=\left\{\gamma\in\mathcal{M}^+(\Omega\times\Omega):\ \gamma(A\times\Omega)=m\mu(A),\ \gamma(\Omega\times B)=m\rho(B) \right\}.
 \end{equation*}
 Of course $\mu\otimes\rho\in\Gamma(\mu,\rho)$. From the classic Optimal Transport theory (\cite{villani,ambrosiogiglisavare,santambrogio-book}) we know that, if $\rho,\mu\in\mathcal{M}^+(\Omega)$ with $m=\rho(\Omega)=\mu(\Omega)$ and $\rho$ is absolutely continuous with respect to the Lebesgue measure, then there is a map $T:\Omega\to\Omega$ such that $T(x)=\nabla\varphi(x)$ for $\rho-$a.e. $x\in\Omega$, $\mu=T\#\rho$ and
 
 \begin{equation}
 \label{eq:wasserstein-metric-classic-1}
 W_{m}^2(\mu,\rho)=\int_{\Omega}|x-\nabla\varphi(x)|^2\rho(x)\ dx.
 \end{equation}
 
  Now, a distance on the space $\Pw{m}{\ }$ is defined as: for $\overrightarrow{\mu},\overrightarrow{\rho}\in\Pw{m}{\ }$ we define the distance by the expression,
 \begin{equation}
 	\label{eq:wasserstein-metric-vector}
 	\dm{m}{\mu}{\rho}:=\sum_{i=1}^N\dmi{m_i}{\mu_i}{\rho_i}.
 \end{equation}
 Note that, due to the compactness of the set $\Omega$, $\Pw{m}{\ }$ is also compact with respect do the weak topology induced by the topological dual of the space $C(\Omega)^N$. This fact implies compactness of $\Pw{m}{\ }$ with respect to the metric $\mathbf{W}_{\overrightarrow{m}}$.
 
 Now, system \eqref{eq:system-origin} has a natural free energy functional associated as in the case of system \eqref{eq:keller-segel}. In order to define this free energy functional take a $\overrightarrow{\rho}\in\mathscr{P}_{\overrightarrow{m}}(\Omega)$ and consider the auxiliary functions $v_i$, $i=1,2,\ldots,N$ defined via,
 \begin{equation}
 \begin{array}{rlr}
 -\Delta v_i(x)=&\rho_i(x),&\text{ for }x\in\text{Int}(\Omega),\\
 v_i(x)=&0,&\text{ for }x\in\partial\Omega.
 \end{array}
 	\label{eq:v_i-auxiliar}
 \end{equation}
In a precise way, $v_i$ is a solution of \eqref{eq:v_i-auxiliar} if the equality
 \begin{align}
 \label{eq:v_i-auxiliar-1}\int_{\Omega}\nabla v_i(x)\cdot\nabla\eta(x)\ dx=&\int_{\Omega}\eta(x)\ d\rho_i(x),
 \end{align}
 holds true for all $\eta\in C^{\infty}_0(\Omega)$. The reader can consult \cite{Augusto-ponce} for the existence and regularity results about equation \eqref{eq:v_i-auxiliar} with a measure data $\rho_i$. With this fact in mind, we define the functional concerned with our system, $\mathcal{F}:\Pw{m}{\ }\to\mathbb{R}\cup\{+\infty\}$, as
 \begin{equation}
 	\label{eq:energy-functional}
 	\mathcal{F}[\overrightarrow{\rho}]=\sum_{i=1}^{N}\mathcal{U}[\rho_i]-\sum_{i,j=1}^{N}\frac{\alpha_i\alpha_j}{2}\int_{\Omega}\nabla v_i(x)\cdot\nabla v_j(x)\ dx,
 \end{equation}
 where $\mathcal{U}$ denotes the Entropy functional,
 \begin{equation*}
 	\label{eq:entropy-functional}
 	\mathcal{U}[\rho]=\int_{\Omega} \rho(x)\log(\rho(x))\ dx,
 \end{equation*}
if $\rho$ is an absolutely continuous positive measure and $\mathcal{U}[\rho]=+\infty$ in another case. It is a well known fact that $\mathcal{F}$ is a Lyapunov functional for the system \eqref{eq:system-origin} without the mass changing term $H$, see \cite{wolansky-multispecies}. In order to use the JKO-scheme, we also consider, for given $\tau>0$ and $\overrightarrow{\nu}\in\Pw{m}{\ }$, the functional $\Phi:\Pw{m}{\ }\to\mathbb{R}\cup\{+\infty\}$ defined via
\begin{equation}
	\label{eq:Moreau-Yosida-functional}
	\Phi(\tau,\overrightarrow{\nu};\overrightarrow{\rho})=\mathcal{F}[\overrightarrow{\rho}]+\frac{1}{2\tau}\dm{m}{\nu}{\rho}.
\end{equation}
In what follows, we use the notation: 
\begin{align}
	\chi=&\left(\max_{1\leq i\leq N}\alpha_i\right)\sum_{i=1}^{N}\alpha_i,\label{eq:notation-chi}\\
	\|\overrightarrow{\rho}\|_{L^{\infty}}=&\max_{i}\|\rho_i\|_{L^{\infty}}.\label{eq:notation-norm-Linfty}
\end{align}
for every $\overrightarrow{\rho}\in\Pw{m}{\ }$ such that $\rho_i\in L^{\infty}(\Omega)$.
\section{$L^{\infty}$ Estimate}
\label{section:Linfty-estimate}
 Following \cite{carrillo-santambrogio}, for given $\overrightarrow{\nu}\in\Pw{m}{\ }$ and $\tau>0$, we consider the optimization problem:

\begin{equation}
	\label{eq:minimization-moreau-yosida}
	 \overrightarrow{\rho} \in\text{Arg}\min\left\{\Phi(\tau,\overrightarrow{\nu};\overrightarrow{\rho}):\overrightarrow{\rho}\in\Pw{m}{\ },\ \|\overrightarrow{\rho}\|_{L^{\infty}}\leq\frac{1}{\chi\tau}\right\}.
\end{equation}
 The next result was shown for the Fokker-Planck equation in \cite[Chapter 8]{santambrogio-book}. We give the proof for the sake of completeness.
 \begin{lemma}
 	\label{lemma:positivity-minimizer}
 	Assume that $\overrightarrow{\rho}$ is a minimizer  in \eqref{eq:minimization-moreau-yosida}. Then, $\rho_i(x)>0$ for almost every where $x\in\Omega$ for all $i=1,\ldots,N$.
 \end{lemma}
 \noindent\emph{Proof.}
 Fix $i_0\in\{1,2,\ldots,N\}$. For $\epsilon\in(0,1)$, define $\overrightarrow{\rho}_{i_0,\epsilon}=(\rho_1,\ldots,\rho_{i_0-1},\rho_{i_0,\epsilon},\rho_{i_0+1},\ldots,\rho_N)$, where $\rho_{i_0,\epsilon}=(1-\epsilon)\rho_{i_0}+\epsilon\tilde{\rho}$, where $\tilde{\rho}=\frac{m_{i_0}}{|\Omega|}1_{\Omega}$ and $|\Omega|$ denotes the Lebesgue measure of the set $\Omega$. Calling $\tilde{v}$ to the solution (in the sense of \eqref{eq:v_i-auxiliar-1}) of $-\Delta\tilde{v}=\tilde{\rho}$ on $\text{Int}(\Omega)$ and $\tilde{v}=0$ on $\partial\Omega$, the minimizing property of $\overrightarrow{\rho}$ can be written as:
 \begin{align*}
 	\mathcal{U}[\rho_{i_0}]-\mathcal{U}[\rho_{i_0,\epsilon}]\leq&\epsilon\alpha_{i_0}\sum_{j=1}^{N}\alpha_j\int_{\Omega}\nabla (v_{i_0}-\tilde{v})\cdot\nabla v_j(x)\ dx+\epsilon^2\alpha_{i_0}^2\int_{\Omega}|\nabla(v_{i_0}-\tilde{v})|^2\ dx\\
 	\ &+\frac{1}{2\tau}\left(\dmi{m_{i_0}}{\nu_{i_0}}{\rho_{i_0,\epsilon}}-\dmi{m_{i_0}}{\nu_{i_0}}{\rho_{i_0}}\right).
 \end{align*}
 Denoting $U(t)=t\log(t)$ that satisfies $U(t)-U(s)\geq U'(s)(t-s)$, for $s>0$, and using the convexity property $W^2_{m_{i_0}}(\nu_{i_0},\rho_{i_0,\epsilon})\leq (1-\epsilon)W^2_{m_{i_0}}(\nu_{i_0},\rho_{i_0})+\epsilon W^2_{m_{i_0}}(\nu_{i_0},\tilde{\rho})$, we arrive to
 \begin{align}
 \int_{\Omega}(\log(\rho_{i_0,\epsilon}(x))+1)(\rho_{i_0}(x)-\tilde{\rho}(x))\ dx \leq&\left[\alpha_{i_0}\sum_{j=1}^{N}\alpha_j\int_{\Omega}\nabla (v_{i_0}-\tilde{v})\cdot\nabla v_j(x)\ dx+\epsilon\alpha_{i_0}^2\int_{\Omega}|\nabla(v_{i_0}-\tilde{v})|^2\ dx\right.\nonumber\\
 \ &\left.+\frac{1}{2\tau}\left(\dmi{m_{i_0}}{\nu_{i_0}}{\tilde{\rho}}-\dmi{m_{i_0}}{\nu_{i_0}}{\rho_{i_0}}\right)\right]\nonumber\\
 =&O(1)\ \ \text{as }\epsilon\to0.\label{eq:lemma-positivity-1}
 \end{align}
 Separating the integral over the sets $A=[\rho_{i_0}>0]$, $B=[\rho_{i_0}=0]$ and noting that
 \begin{equation}
 	(\log(\rho_{i_0,\epsilon}(x))+1)(\rho_{i_0}(x)-\tilde{\rho}(x))\geq (\log(\tilde{\rho}(x))+1)(\rho_{i_0}(x)-\tilde{\rho}(x))\label{eq:lemma-positivity-2},
 \end{equation}
 one gets the inequality,
 \begin{align*}
 -\left(\log\left(\frac{\epsilon m_{i_0}}{|\Omega|}\right)+1\right)\frac{|B|m_{i_0}}{|\Omega|}+\int_{A}(\log(\tilde{\rho}(x))+1)(\rho_{i_0}(x)-\tilde{\rho}(x))\ dx \leq&O(1)\ \ \text{as }\epsilon\to0,
 \end{align*}
 which is a contradiction unless $|B|=0$. It concludes the proof.
 \fin
 Similarly, the next result can be seen as a vector version of \cite[Theorem 1]{carrillo-santambrogio} and we provide the proof as in the previous result. Recall the notation given in \eqref{eq:notation-chi} and\eqref{eq:notation-norm-Linfty}.
 
 \begin{lemma}
 	\label{lemma:existence-solution-Linfty-bound}
 	Problem \eqref{eq:minimization-moreau-yosida} admits a solution. If $\nu_i\in C^{0,a}(\Omega)$, $\inf\nu_i>0$ for all $i=1,\ldots,N$ and $\tau\|\overrightarrow{\nu}\|_{L^{\infty}}\leq4(\lambda-1)/(\chi(2\lambda-1))$ for some $\lambda>1$, then every solution $\overrightarrow{\rho}=(\rho_1,\ldots,\rho_N)$ of \eqref{eq:minimization-moreau-yosida} satisfies that $\rho_i$ is Lipzchits continuous, $\inf\rho_i>0$ for all $i$ and
 	\begin{equation*}
 		\|\overrightarrow{\rho}\|_{L^{\infty}}^{-1}\geq\|\overrightarrow{\nu}\|_{L^{\infty}}^{-1}-\lambda\chi\tau.
 	\end{equation*}
 	Besides that, the identity
 	\begin{equation}
 		\label{eq:euler-lagrange-equation}
 		\nabla\rho_i-\alpha_i\left(\sum_{i=1}^{N}\alpha_j\nabla v_j\right)\rho_i+\frac{\nabla\psi_i}{2\tau}\rho_i=0,
 	\end{equation}
 	holds true a.e. $x\in\Omega$, where $\psi_i$ is the Kantorovich potential between $\rho_i$ and $\nu_i$, this means that $\varphi_i=\frac{|x|^2}{2}-\frac{\psi_i}{2}$ is convex, $\nu_i=\nabla\varphi_i\#\rho_i$ and $\nabla\varphi_i$ satisfies \eqref{eq:wasserstein-metric-classic-1} for $W_{m_i}(\rho_i,\nu_i)$.
 \end{lemma}
 
 \noindent\emph{Proof.}\ \\
 {\it Existence of minimizer:}
 Since the functional $\mathcal{U}$ is classically known to be lower semicontinuous, the details are given only for the second term of the functional $\mathcal{F}$. Fix $i\in\{1,2,\ldots,N\}$. Take $\overrightarrow{\rho}^{k},\overrightarrow{\rho}\in\Pw{m}{\ }\cap(L^{\infty}(\Omega))^N$, such that $\overrightarrow{\rho}^{k}\rightharpoonup\overrightarrow{\rho}$ in the product topology of measures, then it is clear that $\rho_i^k\rightharpoonup\rho_i$ weakly for all $i=1,\ldots, N$. Since, we assume $\|\rho_i^k\|_{L^{\infty}}\leq1/\chi\tau$ then, we can see $(\rho_i^k)_{k\in\mathbb{N}}\subset L^{2}(\Omega)$ as a bounded set. Recalling that $L^{2}(\Omega)\subset\subset H^{-1}(\Omega)$, we can assume, unless subsequence, that $\rho_i^{k}\to f\in H^{-1}(\Omega)$, but by weak convergence necessarily, we have $f=\rho_i$. Thus, the continuity of the operator $(-\Delta)^{-1}:H^{-1}(\Omega)\to H^1(\Omega)$ implies that $v_i^k\to v_i$ in norm, where $-\Delta v_i^k=\rho_i^k$ and $-\Delta v_i=\rho_i$. This fact, obviously, implies
 $$\int_{\Omega}\nabla v_i^k\cdot\nabla v_j^k\ dx\to\int_{\Omega}\nabla v_i\cdot\nabla v_j\ dx, \text{ as }k\to\infty.$$
 Finally, since $\Pw{m}{\ }$ is compact, taking a minimizing sequence it is clear the existence of at least a minimizer $\overrightarrow{\rho}\in\Pw{m}{\ }\cap(L^{\infty}(\Omega))^N$ such that $\|\overrightarrow{\rho}\|_{L^{\infty}}\leq1/\chi\tau$.\\
 
 \noindent{\it Optimality conditions:} As it was conducted in \cite[Chapter 8]{santambrogio-book}, we compute the first variation of $\Phi(\tau,\overrightarrow{\nu};\cdot)$ at the minimizer $\overrightarrow{\rho}$. Fix $i\in\{1,2,\cdots,N\}$. Take $\tilde{\rho}\in\mathcal{M}^+(\Omega)\cap L^{\infty}(\Omega)$ with $\|\tilde{\rho}\|_{L^{\infty}}\leq1/\chi\tau$, $\tilde{\rho}(\Omega)=m_i$ and $\epsilon\in(0,1]$. Considering $\overrightarrow{\rho_{i,\epsilon}}=(\rho_1,\ldots,\rho_{i-1},\rho_{i,\epsilon},\rho_{i+1},\ldots,\rho_N)$, with $\rho_{i,\epsilon}=(1-\epsilon)\rho_i+\epsilon\tilde{\rho}$, we compute the expression:
 $$\lim_{\epsilon\to0}\frac{1}{\epsilon}\left(\Phi(\tau,\overrightarrow{\nu};\overrightarrow{\rho_{i,\epsilon}})-\Phi(\tau,\overrightarrow{\nu};\overrightarrow{\rho_{i}})\right)\geq0.$$
 Considering $\tilde{v}$ such that $-\Delta\tilde{v}=\tilde{\rho}$, we note again that, 
 \begin{align*}
 	\frac{1}{\epsilon}\left(\Phi(\tau,\overrightarrow{\nu};\overrightarrow{\rho_{i,\epsilon}})-\Phi(\tau,\overrightarrow{\nu};\overrightarrow{\rho_{i}})\right)=&\int_{\Omega}\frac{1}{\epsilon}\left(\rho_{i,\epsilon}\log(\rho_{i,\epsilon})-\rho_{i}\log(\rho_{i})\right)\ dx-\alpha_i\sum_{j=1}^{N}\alpha_j\int_{\Omega}\nabla v_j\cdot \nabla(\tilde{v}-v_i)\ dx+O(\epsilon)\\
 	\ &+\frac{1}{2\tau\epsilon}\left(\dmi{m_i}{\nu_{i}}{\rho_{i,\epsilon}}-\dmi{m_i}{\nu_{i}}{\rho_{i}}\right).
 \end{align*}
 To the first integral, we have the bound $\left|\frac{1}{\epsilon}\left(\rho_{i,\epsilon}\log(\rho_{i,\epsilon})-\rho_{i}\log(\rho_{i})\right)\right|\leq\left|\log(\rho_i)+1\right|2/\chi\tau$ and applying Fatou's Lemma to \eqref{eq:lemma-positivity-1}, in view of \eqref{eq:lemma-positivity-2}, it is easy to get that $\log(\rho_i)\in L^1(\Omega)$, so we can use Dominated Convergence theorem. For the last term, Lemma \ref{lemma:positivity-minimizer} gives that $\rho_i>0\ a.e.\ x\in\Omega$ and we can apply \cite[Propositions 7.17 and 7.18]{santambrogio-book} to take the limit as $\epsilon\to0$. In summary, we arrive at
 \begin{equation*}
 	\int_{\Omega}\left(\log(\rho_i)+1-\alpha_i\sum_{j=1}^{N}\alpha_jv_j+\frac{\psi_i}{2\tau}\right)(\tilde{\rho}-\rho_i)\ dx\geq0.
 \end{equation*}
 Taking $l_i:=\text{ess}\inf\left(\log(\rho_i)+1-\alpha_i\sum_{j=1}^{N}\alpha_jv_j+\frac{\psi_i}{2\tau}\right)$ and suitable $\tilde{\rho}$, it follows that
 $$\int_{\Omega}\left(\log(\rho_i)+1-\alpha_i\sum_{j=1}^{N}\alpha_jv_j+\frac{\psi_i}{2\tau}\right)\rho_i\ dx=m_il_i.$$
 So,
 \begin{equation}
 	\label{eq:lemma-existence-solution-Linfty-bound-1}\log(\rho_i)-\alpha_i\sum_{j=1}^{N}\alpha_jv_j+\frac{\psi_i}{2\tau}=l_i-1
 \end{equation}
 for $\rho_i-a.e.\ x\in\Omega$ and, by Lemma \ref{lemma:positivity-minimizer}, for $a.e.\ x\in\Omega$. From regularity theory of the Poisson equation we know that $v_j\in W^{2,p}(\Omega)$, clearly $\psi_i$ is also Lipschitz continuous, then we conclude that $\log(\rho_i)$ is Lipschitz continuous and from an easy computation $\rho_i$ is too. Continuity $\log(\rho_i)$ and compactness of $\Omega$, implies that $\inf\ \rho_i>0$. By deriving \eqref{eq:lemma-existence-solution-Linfty-bound-1} it follows \eqref{eq:euler-lagrange-equation}.
 
 \noindent{\it $L^{\infty} estimate$:} First, by using the regularity result of the Monge-Ampere equation from \cite{caffarelli-mong-amp} we know that $\phi_i=\frac{|x|^2}{2}-\frac{\psi_i}{2}\in C^{2,\bar{a}}$ for $\bar{a}<a$ and so the same is true for the Kantorovich potential $\psi_i$. Denote by $x_i$ the point of minimum of the function $-\alpha_i\sum_{j=1}^{N}\alpha_jv_j+\frac{\psi_i}{2\tau}$ in $\Omega$ which, following \cite{carrillo-santambrogio}, we show that $x_i\in\text{Int}(\Omega)$. Assume $x_i\in\partial\Omega$, condition $\nu_i=\nabla\varphi_i\#\rho_i$, positivity of $\rho_i,\nu_i$ and regularity of $\varphi_i$ implies that $\nabla\varphi_i(\Omega)\subset\Omega$, so $(\nabla\varphi_i(x_i)-x_i)\cdot \overrightarrow{n}(x_i)\leq0$ or $\nabla\psi_i(x_i)\cdot \overrightarrow{n}(x_i)\geq0$. In other hand, since $-\Delta v_j=\rho_j\geq0$ on $\text{Int}(\Omega)$ with $v_j=0$ on $\partial\Omega$, the Minimum principle for sub-harmonic functions implies that all the $v_j'$s reach their minimum at $\partial\Omega$, so by \cite[Lemma 3.4]{gilbarg-truinger} we have
 $$\nabla v_j\cdot\overrightarrow{n}<0, \text{ in }\partial\Omega.$$
 
 Putting together all the previous facts, we get
 $$\left(\frac{1}{2\tau}\nabla\psi_i(x_i)-\alpha_i\sum_{j=1}^{N}\alpha_j\nabla v_j(x_i)\right)\cdot\overrightarrow{n}(x_i)>0,$$
 which contradicts the minimizing property of $x_i$, so $x_i\in\text{Int}(\Omega)$. Note that in this case $x_i$ is a maximum point for $\rho_i$ and therefore $\|\rho_i\|_{L^{\infty}}=\rho_i(x_i)$. It is clear ($v_j\in C^{2,\bar{a}}$) that 
 $$\Delta\psi_i(x_i)\geq2\tau\alpha_i\sum_{j=1}^{N}\alpha_j\Delta v_j(x_i)=-2\tau\alpha_i\sum_{j=1}^{N}\alpha_j\rho_j(x_i).$$
 Recall that the matrix $\nabla^2\phi_i(x_i)=Id_2-\frac{1}{2}\nabla^2\psi_i(x_i)$ is positive. By the inequality $\det(Id_2-\frac{1}{2}\nabla^2\psi_i(x_i))\leq\left(1-\frac{1}{4}\Delta\psi_i(x_i)\right)^2$ and the Monge-Ampere equation one gets,
 \begin{align*}
 	\rho_i(x_i)
 	=&\det\left(Id_2-\frac{1}{2}\nabla^2\psi_i(x_i)\right)\nu_i(\nabla\psi(x_i))\\
 	\leq&\left(1-\frac{1}{4}\Delta\psi_i(x_i)\right)^2\nu_i(\nabla\psi(x_i))\\
 	\leq&\left(1+\frac{1}{2}\tau\alpha_i\sum_{j=1}^{N}\alpha_j\rho_j(x_i)\right)^2\nu_i(\nabla\psi(x_i))\\
 	\leq&\|\overrightarrow{\nu}\|_{L^{\infty}}\left(1+\frac{\tau\chi}{2}\|\overrightarrow{\rho}\|_{L^{\infty}}\right)^2,
 \end{align*}
 or 
 \begin{equation*}
 	\frac{\|\overrightarrow{\rho}\|_{L^{\infty}}}{\left(1+\frac{\tau\chi}{2}\|\overrightarrow{\rho}\|_{L^{\infty}}\right)^2}\leq\|\overrightarrow{\nu}\|_{L^{\infty}}.
 \end{equation*}
 Elementary considerations show that $\displaystyle\frac{s}{1+ks}\leq\frac{s}{(1+s)^2}$ for $s\in[0,k-2]$ where $k>2$. Since the maximum of $\frac{s}{(1+s)^2}$ is reached at $s=1$ and $\tau \chi\|\overrightarrow{\rho}\|_{L^{\infty}}\leq1$, we can conclude that
 \begin{equation*}
 \frac{\|\overrightarrow{\rho}\|_{L^{\infty}}}{1+\frac{\tau\chi k}{2}\|\overrightarrow{\rho}\|_{L^{\infty}}}\leq\|\overrightarrow{\nu}\|_{L^{\infty}},
 \end{equation*}
for  $\tau\|\overrightarrow{\nu}\|_{L^{\infty}}\leq\frac{2(k-2)}{\chi(k-1)^2}$. Taking $\lambda=k/2$ the proof concludes.
 \fin
 
 \begin{remark}
 	In the previous lemma, the restriction $\nu\in C^{0,a}$ can be removed as it was showed in \cite{carrillo-santambrogio}. We prefer to keep this assumption for simplicity.
 \end{remark}
 
 To finish this section, we show a $L^{\infty}$ result for the displacement interpolation
  \begin{lemma}
  	\label{Lemma:Otto}Let $\rho,\nu\in C^{0,\alpha}(\Omega)$ no negative functions with the same mass for some $0<\alpha<1$ and $\inf \rho(x),\inf \nu(x)>0$. Denote for $\varphi$ the convex function such that $\rho=\nabla\varphi\#\nu$ and consider the displacement interpolation
  	$$\rho_{\lambda}=((1-\lambda)x+\lambda\nabla\varphi)\#\nu,$$
  	for $\lambda\in[0,1]$. If $\rho(x),\nu(x)\leq M$ on $\Omega$ for some constant $M>0$ then $\rho_{\lambda}(x)\leq M$ a.e. $x\in\Omega$.
  \end{lemma}
  \noindent\emph{Proof.}
  By invoking regularity theory from \cite{caffarelli-mong-amp} we know that $\varphi\in C^{2,\alpha}(\Omega)$ and the Monge-Ampere equation is satisfied pointwise
  $$\rho(\nabla\varphi(x))\det\nabla^2\varphi(x)=\nu(x).$$
  Thus,
  \begin{align*}
  \rho_{\lambda}((1-\lambda)x+\lambda\nabla\varphi(x))=&\frac{\nu(x)}{\det\left((1-\lambda)I_2+\lambda\nabla^2\varphi\right)}\\
  \leq& \frac{\nu(x)}{\left((1-\lambda)+\lambda(\det\nabla^2\varphi)^{1/2}\right)^2}\\
  =&\frac{\rho(\nabla\varphi(x))\nu(x)}{\left((1-\lambda)\rho(\nabla\varphi(x))^{1/2}+\lambda(\nu(x))^{1/2}\right)^2}\\
  \leq&\frac{\rho(\nabla\varphi(x))\nu(x)}{\left(\rho(\nabla\varphi(x))^{(1-\lambda)/2}(\nu(x))^{\lambda/2}\right)^2}\\
  =&\rho(\nabla\varphi(x))^{\lambda}(\nu(x))^{1-\lambda}\\
  \leq&M.
  \end{align*}
  Since $[(1-\lambda)I-\lambda\nabla\varphi](\Omega)$ is of full measure for $\rho_{\lambda}$ the aim follows.
  \fin
 \section{JKO-Scheme}
 \label{section:JKO-scheme}
 Now, we follow the ideas on \cite{petrelli-tudorasco} in order to write the JKO-scheme in our case. It has been seen that the suitable discrete scheme must be constructed as: given $\overrightarrow{m}^0=(m_1^0,\ldots,m_N^0)\in\mathbb{R}^N$, with $m_i^0>0$ for all $i=1,2,\ldots,N$, and $\overrightarrow{\rho}^0=(\rho_1^0,\ldots,\rho_N^0)\in\Pw{m}{0}$, define
 \begin{align}
 	\overrightarrow{\nu}^k=&\overrightarrow{\rho}^k+\tau\overrightarrow{H}(\overrightarrow{\rho}^k)\label{eq:JKO-1},\\
 	\overrightarrow{m}^{k+1}=&\int_{\Omega}\overrightarrow{\nu}^{k}(x)\ dx\label{eq:JKO-2}\\
 	\overrightarrow{\rho}^{k+1}\in&\text{Arg}\min\left\{\Phi(\tau,\overrightarrow{\nu}^k;\overrightarrow{\rho}):\overrightarrow{\rho}\in\Pw{m}{k+1},\ \|\overrightarrow{\rho}\|_{L^{\infty}}\leq\frac{1}{\chi\tau}\right\}\label{eq:JKO-3},
 \end{align}
for $k=1,2,\ldots$. Now, we make some assumptions about the death and birth rate. This allows to get uniform estimates for the JKO-scheme.
 
 \noindent{\bf H1)}{\it Birth rate:} There is a $C\geq0$ such that $0\leq H_i(\overrightarrow{\rho})\leq C\rho_i$, for $\overrightarrow{\rho}\in\mathbb{R}^N$, $\rho_i\geq0$, $i=1,2,\ldots,N$,\\
 
 \noindent{\bf H2)}{\it Death rate:} There exist constants $0\leq c\leq C\leq1$ such that $-C\rho_i\leq H_i(\overrightarrow{\rho})\leq -c\rho_i$, for $\overrightarrow{\rho}\in\mathbb{R}^N$, $\rho_i\geq0$, $i=1,2,\ldots,N$.
 
 Recall that we have assumed that all the $H_i$ is either a birth rate or death rate. Now, we assert and show the $L^{\infty}$ bound for the discrete solution in the JKO-scheme.
 \begin{theorem}
 	\label{theorem:Linfty-bound-discrete-sol}
 	Let $\overrightarrow{\rho^{_0}}=(\rho_{_1}^{_0},\ldots,\rho_{_N}^{_0})\in\Pw{m}{0}$ and assume that $\rho_i^{_0}\in C^{0,a}(\Omega)$ for some $a\in(0,1)$ and $\inf\rho_i^{_0}>0$ for all $i=1,\ldots,N$. For $\lambda>1$ define $k_0=4(\lambda-1)/(\chi(2\lambda-1))$.
 	\begin{enumerate}
 		\item[a)] {\it Birth rate case:} Assume {\bf H1}, take $T>0$ such that $\chi Te^{CT}<\|\overrightarrow{\rho^{_0}}\|_{L^{\infty}}^{-1}$ and let $\lambda>1$ and $\epsilon_0>0$ be such that $e^{-CT}\|\overrightarrow{\rho^{_0}}\|_{L^{\infty}}^{-1}-\lambda\chi T\geq\epsilon_0$. Take $\tau>0$ such that $\tau(1+C\tau)\leq k_0\epsilon_0$ then, sequence $\overrightarrow{\rho^k}$ in \eqref{eq:JKO-1}-\eqref{eq:JKO-3} is well defined for $k=1,\ldots,\lfloor T/\tau \rfloor$ and satisfies 
 		\begin{equation}
 		\|\overrightarrow{\rho}^k\|_{L^{\infty}}^{-1}\geq (1+C\tau)^{-k}\|\overrightarrow{\rho}^0\|_{L^{\infty}}^{-1}-\lambda\chi\tau k\label{eq:theorem-Linfty-bound-discrete-sol-birth}
 		\end{equation}
 		
 		\item[b)] {\it Death rate case:} Assume {\bf H2}, take $T>0$ such that $\chi T<\|\overrightarrow{\rho}^0\|_{L^{\infty}}^{-1}$ and both $\lambda>1$, $\epsilon_0>0$ satisfies $\|\overrightarrow{\rho}^0\|_{L^{\infty}}^{-1}-\lambda\chi T\geq\epsilon_0$. Take $\tau>0$ such that $\tau\leq k_0\epsilon_0$ then, sequence $\overrightarrow{\rho^k}$ in \eqref{eq:JKO-1}-\eqref{eq:JKO-3} is well defined for $k=1,\ldots,\lfloor T/\tau \rfloor$ and satisfies 
 		\begin{equation}
 		\|\overrightarrow{\rho}^k\|_{L^{\infty}}^{-1}
 		\geq(1-c\tau)^{-k}\left(\|\overrightarrow{\rho}^0\|_{L^{\infty}}^{-1}-\lambda\chi k\tau\right) \label{eq:theorem-Linfty-bound-discrete-sol-death}
 		\end{equation}
 	\end{enumerate}

 \end{theorem}
 \noindent\emph{Proof.} First, we deal with the birth rate case. By induction hypothesis we have that,
 \begin{align*}
 	\|\overrightarrow{\rho}^k\|^{-1}\geq& (1+C\tau)^{-k}\|\overrightarrow{\rho}^0\|^{-1}-\lambda\chi\tau k\\
 	\geq&e^{-CT}\|\overrightarrow{\rho}^0\|^{-1}-\lambda\chi T\geq\epsilon_0.
 \end{align*}
 Therefore, $\tau\|\overrightarrow{\nu}^k\|_{L^{\infty}}\leq\tau(1+C\tau)\|\overrightarrow{\rho}^k\|_{L^{\infty}}\leq\epsilon_0^{-1}\tau(1+C\tau)\leq k_0$. So, we can use Lemma \ref{lemma:existence-solution-Linfty-bound} to ensure the existence of $\overrightarrow{\rho}^{k+1}$. By the same lemma we get,
 \begin{align*}
 	\|\overrightarrow{\rho}^{k+1}\|_{L^{\infty}}^{-1}
 	\geq&(1+C\tau)^{-1}\|\overrightarrow{\rho}^{k}\|_{L^{\infty}}^{-1}-\lambda\chi\tau\\
 	\geq&(1+C\tau)^{-1}\left((1+C\tau)^{-k}\|\overrightarrow{\rho}^0\|^{-1}-\lambda\chi\tau k\right)-\lambda\chi\tau\\
 	\geq&(1+C\tau)^{-(k+1)}\|\overrightarrow{\rho}^0\|^{-1}-\lambda\chi\tau (k+1).
 \end{align*}
 
 Now, for the death rate
 \begin{align*}
 	\|\overrightarrow{\rho}^k\|_{L^{\infty}}^{-1}
 	\geq&(1-c\tau)^{-k}\left(\|\overrightarrow{\rho}^0\|_{L^{\infty}}^{-1}-\lambda\chi k\tau\right)\\
 	\geq&\|\overrightarrow{\rho}^0\|_{L^{\infty}}^{-1}-\lambda\chi T\geq\epsilon_0.
 \end{align*}
 So, we have again, $\tau\|\overrightarrow{\nu}^k\|_{L^{\infty}}\leq\tau(1-c\tau)\|\overrightarrow{\rho}^k\|_{L^{\infty}}\leq\epsilon_0^{-1}\tau\leq k_0$ then, Lemma \ref{lemma:existence-solution-Linfty-bound} applies and
 \begin{align*}
 \|\overrightarrow{\rho}^{k+1}\|_{L^{\infty}}^{-1}
 \geq&(1-c\tau)^{-1}\|\overrightarrow{\rho}^{k}\|_{L^{\infty}}^{-1}-\lambda\chi\tau\\
 \geq&(1-c\tau)^{-(k+1)}\left(\|\overrightarrow{\rho}^0\|_{L^{\infty}}^{-1}-\lambda\chi k\tau\right)-\lambda\chi\tau\\
 \geq&(1-c\tau)^{-(k+1)}\left(\|\overrightarrow{\rho}^0\|_{L^{\infty}}^{-1}-\lambda\chi (k+1)\tau\right),
 \end{align*}
 and we are done.
 \fin
 \begin{remark}
 	\label{remark:Linfty-bound} Under the conditions in Theorem \ref{theorem:Linfty-bound-discrete-sol} we write for the rest of the work that there is $M_T>0$ such that $\|\overrightarrow{\rho^{k}}\|_{L^{\infty}}\leq M_T$ uniformly on $k=1,2,\ldots,\lfloor T/\tau\rfloor$ and $\tau>0$.
 \end{remark}
 \section{Convergence and Weak Solution}
\label{sectio:convergence} 
 In this section, we show the existence of a weak solution of the system \eqref{eq:system-origin} by showing the convergence of the piecewise constant curve given by $\overrightarrow{\rho}^{n}$ defined in \eqref{eq:JKO-1}-\eqref{eq:JKO-3}. In fact, for $\tau>0$ small, we define $\overrightarrow{\rho}_{\tau}:[0,T)\rightarrow (\mathcal{M}(\Omega)\cap L^{\infty}(\Omega))^{N}$ by
 \begin{align}
 	\overrightarrow{\rho}_{\tau}(t)=&\overrightarrow{\rho}^{k+1},\text{ if }t\in[k\tau,(k+1)\tau),\label{eq:aproximate-solution}
 \end{align}
 for $k=1,2,\ldots,\lfloor T/\tau \rfloor$.
 In this direction, one of the most classical estimates, the square sum of the Wasserstein distance, is shown on the next Lemma. Before this, we recall the following result concerned to the Poisson equation \eqref{eq:v_i-auxiliar}, the reader can consult \cite[Proposition 3.2]{Augusto-ponce}.
\begin{proposition}
	\label{prop:existence-posisson-equation}
	Let $\mathcal{O}$ a bounded open set. For every $\rho\in\mathcal{M}(\mathcal{O})$, the linear Dirichlet problem
	$$-\Delta v=\rho,\text{ in } \mathcal{O},$$
	with condition $v=0$ in $\partial\mathcal{O}$, has a solution such that
	$$\|v\|_{L^1}\leq K\|\rho\|_{L^1},$$
	for some constant $K>0$ depending on $\mathcal{O}$.
\end{proposition}
 \begin{lemma}
 	\label{lemma:sum-square-wasserstein}
 	There exists a constant $\bar{C}>0$ that is non depending of $n$ and $\tau>0$, such that
 	$$\sum_{k=0}^{n}\dmk{m}{\nu}{\rho}{k}\leq \bar{C}\tau.$$
 \end{lemma}
 \noindent\emph{Proof.}
 In this result it does not matter to distinguish between the birth/death rate. First, the fact that $\overrightarrow{\rho}^{k+1}$ is a minimizer of $\Phi(\tau,\overrightarrow{\nu}^{k};\cdot)$ implies that,
 \begin{align*}
 	\frac{1}{2\tau}\dmk{m}{\nu}{\rho}{k}
 	\leq& \mathcal{F}(\overrightarrow{\nu^{k}})-\mathcal{F}(\overrightarrow{\rho^{k+1}})\\
 	=&\sum_{i=1}^{N}\left(\mathcal{U}(\nu_i^{k})-\mathcal{U}(\rho_i^{k+1})\right)-\sum_{i,j=1}^{N}\alpha_i\alpha_j\left(\int_{\Omega}\nabla \tilde{v}_i^{k}\cdot\nabla \tilde{v}_j^{k}-\nabla v_i^{k+1}\cdot\nabla v_j^{k+1}\right)\ dx\\
 	=:&A_k+B_k,
 \end{align*}
 where $\tilde{v}_i^{k}=v_i^{k}+\tau g_i^{k}$, $-\Delta g_i^{k}=H_i(\overrightarrow{\rho}^k)$ on $\text{Int}(\Omega)$ and $g_i^{k}=0$ on $\partial\Omega$. We estimate $A_k$, by noting that $\nu_i^{k}=(1-\tau)\rho_i^k+\tau(\rho_i^k+H_i(\overrightarrow{\rho}^k))$. Thus,
 \begin{align*}
 	\mathcal{U}(\nu_i^{k})-\mathcal{U}(\rho_i^{k+1})\leq&\tau\left(\mathcal{U}(\rho_i^k+H_i({\overrightarrow{\rho}^k}))-\mathcal{U}(\rho_i^k)\right)+\mathcal{U}(\rho_i^{k})-\mathcal{U}(\rho_i^{k+1}).
 \end{align*}
 Now from Theorem \ref{theorem:Linfty-bound-discrete-sol} (see Remark \ref{remark:Linfty-bound}), we know that $\rho_i^k(x)\in[0,M_T]$, so we have $\rho_i^k(x)+H_i(\overrightarrow{\rho}^k(x))\in I_T:=[0,(1+C)M_T]$. Recalling that $U(s)=s\log(s)$ we arrive at,
 \begin{equation}
 	\sum_{k=0}^{n}A_k\leq\tau(n+1) N|\Omega|(\sup_{I_T}U-e^{-1})+\sum_{i=1}^{N}\mathcal{U}(\rho_i^0)-N|\Omega|e^{-1}.\label{eq:lemm-sum-square-wasserstein-1}
 \end{equation}
 
 For $B_k$, we start by noting
 \begin{align*}
 	\int_{\Omega}\nabla\tilde{v}_i^{k}\cdot\nabla\tilde{v}_j^{k}\ dx
 	=&\int_{\Omega}\nabla v_i^{k}\cdot\nabla v_j^{k}\ dx+\tau\left(\int_{\Omega}\nabla v_i^{k}\cdot\nabla g_j^{k}\ dx+\int_{\Omega}\nabla g_i^{k}\cdot\nabla v_j^{k}\ dx\right)\\
 	\ &+\tau^2\int_{\Omega}\nabla g_i^{k}\cdot\nabla g_j^{k}\ dx.
 \end{align*}
 Then, we resort to Proposition \ref{prop:existence-posisson-equation} and it follows,
 \begin{align*}
 	-\int_{\Omega} \nabla v_i^k\cdot\nabla g_j^k\ dx
 	\leq&C\|v^k_i\|_{L^1}\|\rho_j^k\|_{L^{\infty}}\\
 	\leq&CKm_i^k\|\rho_j^k\|_{L^{\infty}},
 \end{align*}
 similarly,
 \begin{align*}
 -\int_{\Omega} \nabla g_i^k\cdot\nabla g_j^k\ dx
 \leq&C\|g^k_i\|_{L^1}\|\rho_j^k\|_{L^{\infty}}\\
 \leq&C^2Km_i^k\|\rho_j^k\|_{L^{\infty}}.
 \end{align*}
 
 Join all these estimates, it follows that
 \begin{align*}
 	B_k\leq&-\sum_{i,j=1}^{N}\alpha_i\alpha_j\int_{\Omega}\left(\nabla v_i^{k}\cdot\nabla v_j^{k}-\nabla v_i^{k+1}\cdot\nabla v_j^{k+1}\right)\ dx+(2\tau CK\chi+\tau^2C^2K\chi)\left(\sum_{i=1}^Nm_i^{k}\right)\|\overrightarrow{\rho}^k\|_{L^{\infty}}\\
 	\leq&-\sum_{i,j=1}^{N}\alpha_i\alpha_j\int_{\Omega}\left(\nabla v_i^{k}\cdot\nabla v_j^{k}-\nabla v_i^{k+1}\cdot\nabla v_j^{k+1}\right)\ dx+\tau CK\chi\left(2+\tau C\right)\left(\sum_{i=1}^Nm_i^{0}\right)(1+C\tau)^{k}\|\overrightarrow{\rho}^k\|_{L^{\infty}}\\
 	\leq&-\sum_{i,j=1}^{N}\alpha_i\alpha_j\int_{\Omega}\left(\nabla v_i^{k}\cdot\nabla v_j^{k}-\nabla v_i^{k+1}\cdot\nabla v_j^{k+1}\right)\ dx+\tau CK\chi\left(2+ C\right)\boldsymbol{m}^{0}e^{Ct}M_T,
 \end{align*}
 where $k\tau\leq t<(k+1)\tau$ and $\boldsymbol{m}^{0}=\sum_{i=1}^Nm_i^{0}$. This last inequality implies the one
 \begin{align}
 	\sum_{k=0}^{n}B_k
 	\leq&-\sum_{i,j=1}^{N}\alpha_i\alpha_j\int_{\Omega}\left(\nabla v_i^{0}\cdot\nabla v_j^{0}-\nabla v_i^{n+1}\cdot\nabla v_j^{n+1}\right)\ dx+\tau(n+1) CK\chi\left(2+ C\right)\boldsymbol{m}^{0}e^{Ct}M_T.\label{eq:lemm-sum-square-wasserstein-2}
 \end{align}
 The same arguments imply the estimate
 \begin{align*}
 	\sum_{i,j=1}^{N}\alpha_i,\alpha_j\int_{\Omega}\nabla v_i^{n+1}\cdot\nabla v_j^{n+1}\ dx\leq&K\chi(1+C)e^{Ct}\boldsymbol{m}^{0}M_T.
 \end{align*}
 The proof concludes taking into account \eqref{eq:lemm-sum-square-wasserstein-1}, \eqref{eq:lemm-sum-square-wasserstein-2} and the previous inequality.
 \fin
 Now we are going to show the convergence of the curve $\overrightarrow{\rho}_{\tau}$ defined in \eqref{eq:aproximate-solution} and some kind of regularity of its limit. Moreover, we also show that this curve is a weak solution of the system \eqref{eq:system-origin}. This proof follows the main ideas of the JKO scheme (\cite{JKO,ambrosiogiglisavare}) but with the technical adaptations that are required.

 \begin{theorem}
 	\label{theorem:convergence} Under the conditions in Theorem \ref{theorem:Linfty-bound-discrete-sol}, we have $\overrightarrow{\rho}_{\tau}(t,\cdot)\in (W^{1,p}(\Omega))^N$ for all $t\in(0,T]$ and $1\leq p\leq\infty$.  Moreover, for $p=2$, it follows that $\nabla\overrightarrow{\rho}_{\tau}(\cdot)\in (L^2([0,T]\times\Omega))^{2N}$ and $\overrightarrow{\rho}_{\tau}(\cdot)\in (L^2([0,T]\times\Omega))^N$ are weakly compact as $\tau\to0$ and for any limits $\overrightarrow{\rho}(\cdot)$ of $\overrightarrow{\rho}_{\tau}(\cdot)$ and $G(\cdot)$ of $\nabla\overrightarrow{\rho}_{\tau}(\cdot)$ we have $G(t,x)=\nabla\overrightarrow{\rho}(t,x)$.
 \end{theorem}
 \noindent\emph{Proof:} We separate the proof in some steps:\\
 {\it $L^p$-estimates:} By the regularity stated in Lemma \ref{lemma:existence-solution-Linfty-bound} we can differentiate the optimality conditions \eqref{eq:lemma-existence-solution-Linfty-bound-1} and multiply by $\rho_i^{k+1}$. So, for $\xi\in C^{\infty}_c(\Omega;\mathbb{R}^2)$ we get
 \begin{align}
 	\int_{\Omega}\nabla\rho_i^{k+1}\cdot\xi\ dx=&-\int_{\Omega}\frac{x-\nabla\phi_i(x)}{\tau}\cdot\xi\rho_i^{k+1}\ dx+\alpha_i\int\sum_{j=1}^{N}\alpha_j\nabla v_j^{k+1}\cdot \xi\ \rho_i^{k+1}\ dx\label{eq:theorem-convergence-1}\\
 	\leq&\frac{\text{diam}(\Omega)}{\tau}\|\rho_i^{k+1}\|_{L^{p}}\|\xi\|_{L^{p'}}+\alpha_i\sum_{j=1}^{N}\alpha_j\left\|\nabla v_j^{k+1}\right\|_{L^{p}}\|\xi\|_{L^{p'}}M_T\nonumber\\
 	\leq&\left(\frac{\text{diam}(\Omega)}{\tau}+\alpha_iC\sum_{j=1}^{N}\alpha_jM_T\right)\|\rho_i^{k+1}\|_{L^{p}}\|\xi\|_{L^{p'}},\nonumber
 \end{align}
 where in the last inequality we have used Calderón-Zygmund theory. So, $\rho_{i}^{k+1}\in W^{1,p}(\Omega)$ and 
 \begin{equation*}
 	\|\nabla\rho_{i}^{k+1}\|_{L^{p}}\leq \left(\frac{\text{diam}(\Omega)}{\tau}+C\chi M_T\right)\|\rho_{i}^{k+1}\|_{L^{p}},
 \end{equation*}
 for all $1\leq p\leq\infty$. In the case, when $p=2$ the estimate can be refined by noting,
 \begin{align*}
 	-\int_{\Omega}\frac{x-\nabla\phi_i(x)}{\tau}\cdot\xi(x)\rho_i^{k+1}\ dx\leq&\left(\int_{\Omega}\frac{|x-\nabla\phi_i(x)|^2}{\tau^2}\rho_i^{k+1}\ dx\right)^{1/2}\left(\int_{\Omega}|\xi(x)|^2\rho_i^{k+1}\ dx\right)^{1/2}\\
 	\leq&\frac{1}{\tau}W_{m_i}(\nu_i^{k},\rho_i^{k+1})\|\rho_i^{k+1}\|_{L^{\infty}}^{1/2}\|\xi\|_{L^{2}}\\
 	\leq&\frac{M_T^{1/2}}{\tau}W_{m_i}(\nu_i^{k},\rho_i^{k+1})\|\xi\|_{L^{2}}.
 \end{align*}
 So, by inserting this last inequality into \eqref{eq:theorem-convergence-1}, we have
 \begin{align*}
 	\|\nabla\rho_{i}^{k+1}\|_{L^{2}}\leq \frac{M_T^{1/2}}{\tau}W_{m_i}(\nu_i^{k},\rho_i^{k+1})+C\chi M_T\|\rho_{i}^{k+1}\|_{L^{2}}.
 \end{align*}
 This last estimate, allows to get the next one
 \begin{equation}
 	\label{eq:theorem-convergence-L2-gradient}
 	\|\nabla\rho_{i}^{k+1}\|_{L^{2}}^2\leq C_T\left(\frac{W^2_{m_i}(\nu_i^{k},\rho_i^{k+1})}{\tau^2}+1\right),
 \end{equation}
 for some $C_T>0$.
 
 \noindent{\it Weak $L^2$ compactness:} Denote $Q=[0,T]\times \Omega$. The previous step and Lemma \ref{lemma:sum-square-wasserstein} implies that,
 \begin{equation*}
 	\int\!\int_{Q}|\nabla\rho_{i,\tau}(t,x)|^2\ dx\ dt\leq C_T(C+T),
 \end{equation*}
 for all $\tau>0$ and it is clear from Theorem \ref{theorem:Linfty-bound-discrete-sol} that $\rho_{i,\tau}(t,x)\in L^{2}(Q)$ is also bounded uniformly with respect to $\tau$. Therefore, there are functions $\rho_i\in L^2(Q)$ and $G_i\in (L^{2}(Q))^2$ and a sequence $\tau_l\to0$ such that $\rho_{i,\tau_l}\to \rho_i$, $\nabla\rho_{i,\tau_l}\to G_i$ weakly in $L^2$. But by taking test functions $\xi\in C^{\infty}_c(\text{Int}(Q))$ and using the $L^2$ weak convergence, it follows that $G_i=\nabla\rho_i$.
 \fin
 
 In view of the non linearity on system \eqref{eq:system-origin}, we need a stronger compactness result in order to take the limit to the discrete solution. This is the reason for which we adapt the arguments of \cite{petrelli-tudorasco,Otto} to our setting.
 
 \begin{theorem}
 	\label{theorem:strong-compact} The family $\overrightarrow{\rho}_{\tau}(\cdot)\in (L^1([0,T]\times\Omega))^N$, $\tau>0$ is strongly compact.
 \end{theorem}
 \noindent\emph{Proof.}
 The idea is to apply the Riesz-Fréchet-Kolmogorov criterion for strong compactness. Define $P_{\tau,i}:=\rho_{\tau,i}$ on $[0,T]\times\Omega$ and $P_{\tau,i}:=0$ outside the $\mathbb{R}^3\setminus[0,T]\times\Omega$. Let $s\in\mathbb{R}$ and $h\in\mathbb{R}^2$, there is no matter in assuming that $s>0$. So,
 \begin{align*}
 	\int\!\!\!\int_{[0,T]\times\Omega} |P_{\tau,i}(t+s,x+h)-P_{\tau,i}(t,x)|\ dxdt
 	\leq& \int_s^{T}\!\!\!\int_{\Omega} |P_{\tau,i}(t,x+h)-P_{\tau,i}(t,x)|\ dxdt\\
 	\ &+\int_{0}^{T-s}\!\!\!\int_{\Omega} |\rho_{\tau,i}(t+s,x)-\rho_{\tau,i}(t,x)|\ dxdt\\
 	\ &+\int_{T-s}^{T}\!\!\!\int_{\Omega} \rho_{\tau,i}(t,x)\ dxdt\\
 	=&:I_1+I_2+I_3.
 \end{align*}
 
 \noindent {\it Estimate for $I_1$}: Since $\rho_{\tau,i}$ is Lipschitz in the spatial variable and using the convexity of $\Omega$ we arrive at
 \begin{align*}
 	I_1=&\int_s^{T}\!\!\!\int_{\Omega\cap(\Omega-h)} |\rho_{\tau,i}(t,x+h)-\rho_{\tau,i}(t,x)|\ dxdt\\
 	\ &+\int_s^{T}\!\!\!\int_{\Omega\cap(\Omega-h)^c} \rho_{\tau,i}(t,x)\ dxdt\\
 	\leq&|h|\int_s^{T}\!\!\!\int_{\Omega\cap(\Omega-h)}\!\!\int_{0}^{1} |\nabla\rho_{\tau,i}(t,x+\lambda h)|\ d\lambda dxdt\\
 	+\ &TM_T|\Omega\cap(\Omega-h)^c|\\
 	=&|h|\int_s^{T}\!\!\!\int_{0}^{1}\!\!\int_{\Omega\cap(\Omega-h)+\lambda h} |\nabla\rho_{\tau,i}(t,y)|\ dyd\lambda dt\\
 	+\ &TM_T|\Omega\cap(\Omega-h)^c|\\
 	\leq&|h|\int_0^{T}\!\!\!\int_{\Omega} |\nabla\rho_{\tau,i}(t,y)|\ dy dt\\
 	+\ &TM_T|\Omega\cap(\Omega-h)^c|\to0\\
 \end{align*}
 as $|h|\to0$ uniformly with respect to $\tau$.
 
 \noindent{\it Estimate for $I_2$}: Since we need a uniform estimate with respect to $\tau$, it is convenient to take $j\geq0$ an integer such that $s=j\tau+r$ with $0\leq r<\tau$. Note that $r$ is depending on $\tau$ and $s$. Assuming that $n\tau\leq T <(n+1)\tau$ we can compute
 \begin{align}
 	I_2\leq&\sum_{k=0}^{n-j}\int_{k\tau}^{(k+1)\tau}\!\!\!\int_{\Omega}|\rho_{\tau,i}(t+s,x)-\rho_{\tau,i}(t,x)|^2\ dxdt\nonumber\\
 	\leq&(\tau-r)\sum_{k=0}^{n-j}\int_{\Omega}|\rho_{i}^{k+j+1}(x)-\rho_{i}^{k+1}(x)|^2\ dx+r\sum_{k=0}^{n-j}\int_{\Omega}|\rho_{i}^{k+j+2}(x)-\rho_{i}^{k+1}(x)|^2\ dx.\label{eq:theorem-strong-compact-1}
 \end{align}
 In order to estimate \eqref{eq:theorem-strong-compact-1}, we call $\xi(x)=\rho_{i}^{k+j+1}(x)-\rho_{i}^{k+1}(x)$ and write
 \begin{align}
 	|\rho_{i}^{k+j+1}(x)-\rho_{i}^{k+1}(x)|^2=&\sum_{l=1}^{j}(\rho_{i}^{k+l+1}(x)-\nu_{i}^{k+l}(x))\xi(x)+\sum_{l=1}^{j}(\nu_{i}^{k+l}(x)-\rho_{i}^{k+l}(x))\xi(x)\nonumber\\
 	=&\sum_{l=1}^{j}(\rho_{i}^{k+l+1}(x)-\nu_{i}^{k+l}(x))\xi(x)+\sum_{l=1}^{j}\tau H_i(\overrightarrow{\rho}^{k+l})\xi(x).\label{eq:theorem-strong-compact-2}
 \end{align}
 Then the first term is treated as:
 \begin{align*}
 	\int_{\Omega}(\rho_{i}^{k+l+1}(x)-\nu_{i}^{k+l}(x))\xi(x)\ dx=&\int_{\Omega}(\xi(x)-\xi(\nabla\varphi_{i}^{k+l}(x)))\rho_{i}^{k+l+1}(x)\ dx\\
 	\leq&\int_{0}^{1}\!\!\!\int_{\Omega}|\nabla\xi((1-\lambda)x+\lambda\nabla\varphi_{i}^{k+l}(x))||x-\nabla\varphi_{i}^{k+l}(x)|\rho_{i}^{k+l+1}(x)\ dx\ d\lambda\\
 	\leq&\left(\int_{0}^{1}\!\!\!\int_{\Omega}|\nabla\xi((1-\lambda)x+\lambda\nabla\varphi_{i}^{k+l}(x))|^2\rho_{i}^{k+l+1}(x)\ dx\ d\lambda\right)^{1/2} W_{m_i}(\rho_i^{k+l+1},\nu_i^{k+l})\\
 	\leq&M_T^{1/2}\|\nabla(\rho^{k+j+1}(x)-\rho^{k+1}(x))\|_{L^2}W_{m_i}(\rho_i^{k+l+1},\nu_i^{k+l}),
  \end{align*}
  where in the last inequality we have used Lemma \ref{Lemma:Otto}. Thus, by using Theorem \ref{theorem:convergence} the last estimate continues as,
  \begin{align*}
  	\int_{\Omega}(\rho_{i}^{k+l+1}(x)-\nu_{i}^{k+l}(x))\xi(x)\ dx
  	\leq&M_T\frac{W_{m_i}(\rho_i^{k+l+1},\nu_i^{k+l})}{\tau}\left(W_{m_i}(\rho_i^{k+j+1},\nu_i^{k+j})+W_{m_i}(\rho_i^{k+1},\nu_i^{k})+2\tau C\right)\\
  	\leq& \frac{C}{2\tau}\left(2W_{m_i}^2(\rho_i^{k+l+1},\nu_i^{k+l})+W_{m_i}^2(\rho_i^{k+j+1},\nu_i^{k+j})\right.\\
  	\ &\left.+W_{m_i}^2(\rho_i^{k+1},\nu_i^{k})+2\tau W_{m_i}(\rho_i^{k+l+1},\nu_i^{k+l})\right).
  \end{align*}
  Inserting the last estimate and \eqref{eq:theorem-strong-compact-2} into \eqref{eq:theorem-strong-compact-1} and using Lemma \ref{lemma:sum-square-wasserstein}, we arrive at
  \begin{align*}
  	I_2\leq&(\tau-r)\sum_{l=1}^{j}C(1+\sqrt{\tau})+rC\sum_{l=1}^{j}(1+\sqrt{\tau})+C\tau^2j(n-j)\\
  	\leq&sC.
  \end{align*}
  For a constant $C>0$ that does not depend on $s$ and $\tau$.\\
  \noindent{\it Estimate for $I_3$}: It is an easy consequence of the $L^{\infty}$ estimate of $\rho_{\tau,i}$, that
  
 $$I_3\leq Cs.$$
The proof concludes for applying the Riesz-Fréchet-Kolmogorov criterion, see \cite[Theorem 4.26]{brezis}.
 \fin
 
 Now, we are going to show that the limit $\overrightarrow{\rho}=(\rho_1,\ldots,\rho_N)$ obtained in Theorem \ref{theorem:convergence} is a weak solution of \eqref{eq:system-origin}.
 \begin{theorem}
 	\label{theorem:weak-solution}
 	Let $\overrightarrow{\rho}(t)=(\rho_1(t),\ldots,\rho_N(t))$ the curve given by Theorem \ref{theorem:convergence} and assume that $H$ is Lipschitz continuous. Then, $\overrightarrow{\rho}(t)$ is a weak solution of the problem \eqref{eq:system-origin}.
 \end{theorem}
 
 \noindent\emph{Proof.}
 We follow the standard arguments for the step descendant scheme. It follows, by multiplying \eqref{eq:euler-lagrange-equation} by $\nabla\xi$ for a test function $\xi\in C^{\infty}_c(\text{Int}(Q))$ and integrating over $\Omega$, that
 \begin{equation*}
 	\int_{\Omega}\left(\nabla\rho_i^{k+1}\cdot\nabla\xi(t,x)-\alpha_i\left(\sum_{i=1}^{N}\alpha_j\nabla v_j^{k+1}\cdot\nabla\xi(t,x)\right)\rho_i^{k+1}\right)\ dx+\int_{\Omega}\frac{(x-\nabla\varphi_i^{k+1}(x))}{\tau}\cdot\nabla\xi(t,x)\rho_i^{k+1}\ dx=0.
 \end{equation*}
 Taylor expansion $\xi(t,\nabla\varphi_i^{k+1}(x))-\xi(t,x)=-\nabla\xi(t,x)\cdot(x-\nabla\varphi_i^{k+1}(x))+O(|x-\nabla\varphi_i^{k+1}(x)|^2)$ and the identity $\nabla\varphi_i^{k+1}\#\nu_i^k=\rho_i^{k+1}$ leads 
 \begin{align*}
 	\int_{\Omega}\frac{(x-\nabla\varphi_i^{k+1}(x))}{\tau}\cdot\nabla\xi(x)\rho_i^{k+1}\ dx
 	=&\frac{1}{\tau}\left(\int_{\Omega}\xi(t,x)\rho_i^{k+1}(x)\ dx-\int_{\Omega}\xi(t,x)\rho_i^{k}(x)\ dx\right)\\
 	&-\int_{\Omega}\xi(t,x)H_i(\overrightarrow{\rho}^{k+1}(x))\ dx+O(\dmi{m_i}{\nu_i^{k}}{\rho_i^{k+1}}/\tau).
 \end{align*}
 Integrating over $[k\tau,(k+1)\tau)$ and adding from $k=1$ to $n$, one get
 \begin{align}
 	\int_{\tau}^{n\tau}&\int_{\Omega}\left(\nabla\rho_{i,\tau}(t,x)\cdot\nabla\xi(t,x)-\alpha_i\sum_{i=1}^{N}\alpha_j\nabla v_{j,\tau}(t,x)\cdot\nabla\xi(t,x)\rho_{i,\tau}(t,x)\right)\ dx\ dt\nonumber\\
 	\ &\int_{\tau}^{n\tau}\int_{\Omega}\rho_{i,\tau}(t,x)\frac{\xi(t+\tau,x)-\xi(t,x)}{\tau}\ dx\ dt+\int_{\tau}^{n\tau}\xi(t,x)H_i(\overrightarrow{\rho}_{\tau}(t,x))\ dx\ dt=O(\tau),\label{eq:weak-solution}
 \end{align}
 for $\tau>0$ small enough. Note the facts: $\rho_{i,\tau}\to\rho_{i}$ $L^1$ strongly, $\nabla\rho_{i,\tau}\to\nabla\rho_{i}$ weakly in $L^2$ and by Proposition \ref{prop:existence-posisson-equation} it follows that $v_{j,\tau}\to v_j$ strongly in $L^1(Q)$ for some $v_j(\cdot,t)\in L^1(\Omega)$ and the $L^2$ compactness in the spatial variables of $\rho_{i,\tau}$ got in Theorem \ref{theorem:strong-compact} implies that $-\Delta v_j=\rho_j$ in $\text{Int}(\Omega)$ and $v_j=0$ on $\partial\Omega$. Putting all these facts when we take the limit  as $\tau\to0$ in \eqref{eq:weak-solution} for some subsequence, it gives,
  \begin{align*}
  \int_{0}^{T}\int_{\Omega}\left(\nabla\rho_{i}(x,t)\cdot\nabla\xi(x,t)-\alpha_i\sum_{i=1}^{N}\alpha_j\nabla v_{j}(x,t)\cdot\nabla\xi(x,t)\rho_{i,\tau}(x,t)+\rho_{i}(x,t)\frac{\partial\xi(x,t)}{\partial t}\right)\ dx\ dt\\
 +\int_{0}^{T}\int_{\Omega}\xi(x,t)H_i(\overrightarrow{\rho}(x,t))\ dx\ dt=0,
  \end{align*}
 that concludes the proof.
 \fin
 
 \section{Blowing-up of Solutions}
\label{section:blow-up} 
 In this last section, we study the phenomenon of mass concentration. It was already prevent by Theorem \ref{theorem:Linfty-bound-discrete-sol} that $L^{\infty}$ solutions could not exist globally in time. Indeed, we show that for radially symmetric solutions blowing-up occurs if some conditions are given on the parameters. We concentrate on the case with two species and the domain being the unitary disk $\Omega=\bar{D}:=\overline{D(0,1)}$:
 \begin{equation}
 \left\{
 \begin{array}{rcll}
 \displaystyle\frac{\partial \rho_i}{\partial t}&=&\displaystyle\Delta\rho_i-\alpha_i\nabla \cdot\left(\rho_i\nabla v\right)-c_i\rho_i,&\text{ for }(x,t)\in D(0,1)\times (0,T),\\
  -\Delta v&=&\alpha_1\rho_1+\alpha_2\rho_2,&\text{ for }(x,t)\in D(0,1)\times (0,T),\\
 (\alpha_i\rho_i\nabla v-\nabla\rho_i)\cdot \overrightarrow{n}&=&0,\ v=0&\text{ for } (x,t)\in \partial D(0,1)\times (0,T),\\
 \rho_i(x,0)&=&\rho_{i,0}(x)\geq0&\text{ for }x\in D(0,1),
 \end{array}
 \right.\label{eq:system-origin-two-species}
 \end{equation}
 for $i=1,2$ and $c_1,c_2>0$ being constants. Now, define the quantities
 \begin{align}
 	M_i(t)=e^{c_it}\int_{\bar{D}} |x|^2\rho_i(x,t)\ dx\ \text{ for }i=1,2.
 \end{align}
 First, we show the following technical lemma in order to analyzing the blowing-up of the solutions.
 \begin{lemma}
 	\label{lemma:estimate-second-moment}
 	Let us consider a radial classic solution $\rho_i:\bar{D}\times[0,T)\to\mathbb{R}$, $i=1,2$ for the system \eqref{eq:system-origin-two-species}. Then, we have the estimate
 	\begin{equation}
 		\label{eq:lemma-estimate-second-moment}
 		e^{-c_1t}M_1'+e^{-c_2t}M_2'\leq 4\left(m_{1,0}e^{-c_1t}+m_{2,0}e^{-c_2t}\right)-\frac{1}{2\pi}\left(\alpha_1m_{1,0}e^{-c_1t}+\alpha_2m_{2,0}e^{-c_2t}\right)^2,
 	\end{equation}
 	where $m_{i,0}>0$ is the initial mass for $\rho_i$ and $c_i\geq0$, $i=1,2$ are constants.
 \end{lemma}
\noindent\emph{Proof:}
 	The proof is actually standard and we give this for the sake of completeness. It is easy to see that
 	$$\int_{\bar{D}}\rho_i(x,t)\ dx=m_{i,0}e^{-c_it}.$$
 	Now by multiplying \eqref{eq:system-origin-two-species} by $|x|^2$ and using by parts integration, we get
 	\begin{align*}
 		\frac{d}{dt}\int_{\bar{D}}|x|^2\rho_i(x,t)\ dx\leq&\ 4m_{i,0}e^{-c_it}+2\alpha_i\int_{\bar{D}}\rho_{i}(x,t)\nabla v\cdot x\  dx-c_i\int_{\bar{D}}|x|^2\rho_i(x,t)\ dx
 	\end{align*}
 	or
 	\begin{align*}
 	M_i'(t)\leq&\ 4m_{i,0}+2e^{c_it}\alpha_i\int_{\bar{D}}\rho_{i}(x,t)\nabla v\cdot x\  dx.
 	\end{align*}
	 Now, due that solutions are radial, it follows that $x\cdot \nabla v=r\frac{\partial v}{\partial r}$ and $-\frac{1}{r}\frac{\partial}{\partial r}\left(r\frac{\partial v}{\partial r}\right)=\alpha_1\rho_1(r,t)+\alpha_2\rho_2(r,t)$. Thus, we can estimate,
	 \begin{align*}
	 M_i'(t)\leq&\ 4m_{i,0}-4\pi e^{c_it}\alpha_i\int_{0}^{1}\rho_{i}(r,t)\left[\int_{0}^{r}s\left(\alpha_1\rho_1(s,t)+\alpha_2\rho_2(s,t)\right)\ ds\right] r\ dr.
	 \end{align*}
	 Multiplying by $e^{-c_it}$ and adding over $i=1,2$ the inequality follows.
 \fin
We adopt the following notation
\begin{equation}
\label{eq:critical-constant}
\Lambda=\frac{(\alpha_1m_{1,0}+\alpha_2m_{2,0})^2}{8\pi(m_{1,0}+m_{2,0})}.
\end{equation}
 
 \begin{theorem}
 	\label{theorem:blow-up-solutions-1}
 	Let us denote by $\rho_i(x,t)$ a radially symmetric and classic solution of \eqref{eq:system-origin-two-species}. Denote $C=\max\{c_1,c_2\}$ and $c=\min\{c_1,c_2\}$ and assume that $c_1,c_2>0$,
 	\begin{equation}
 	\label{eq:critical-inequality}
 		\Lambda>1\text{ and } \frac{2c}{C}-4\leq\Lambda^{\frac{c}{2C-c}}\left(\frac{2c}{C}\Lambda-4-c\right).
 	\end{equation}
 	Then the maximal time of existence $T_{\max}$ is finite and solutions present blowup in finite time.
 \end{theorem}
 
 \noindent\emph{Proof:}
We suppose that solution is global in time and then we apply Lemma \ref{lemma:estimate-second-moment}. So, we get
 \begin{align}
 	e^{-c_1t}M_1'+e^{-c_2t}M_2'\leq 4\left(m_{1,0}+m_{2,0}\right)e^{-ct}-\frac{1}{2\pi}\left(\alpha_1m_{1,0}+\alpha_2m_{2,0}\right)^2e^{-2Ct}.\label{eq:theorem-blow-up-soutions-1-1}
 \end{align}
 Now, we use the following notation 
 $$\beta(t)=\int_{\bar{D}}|x|^2\left(\rho_1(x,t)+\rho_2(x,t)\right)\ dx+\int_{0}^{t}\int_{\bar{D}}|x|^2\left(c_1\rho_1(x,\tau)+c_2\rho_2(x,\tau)\right)\ dx\ d\tau.$$
 Integrating \eqref{eq:theorem-blow-up-soutions-1-1} from $0$ to $t$ we arrive to,
 \begin{align*}
 	\beta(t)\leq& \left[\frac{e^{-2Ct}}{4\pi C}\left(\alpha_1m_{1,0}+\alpha_2m_{2,0}\right)^2-\frac{4e^{-ct}}{c}\left(m_{1,0}+m_{2,0}\right)\right]-\left[\frac{1}{4\pi C}\left(\alpha_1m_{1,0}+\alpha_2m_{2,0}\right)^2-\frac{4}{c}\left(m_{1,0}+m_{2,0}\right)\right]\\
 	\ &+m_{1,0}+m_{2,0}.
 \end{align*}
 
 We are going to compute the minimum of the right hand side of the last inequality. In fact, it is easy to see that the point of the minimum $t_{min}$ satisfies the relation:
 \begin{equation*}
 	e^{(2C-c)t_{min}}=\Lambda.
 \end{equation*}
 Thus,
 \begin{align*}
 	\beta(t_{min})\leq& \left[e^{-ct_{min}}\left(\frac{2}{C}-\frac{4}{c}\right)-\frac{2}{C}\Lambda+\frac{4}{c}-1\right]\left(m_{1,0}+m_{2,0}\right)\\
 	=&\left[\Lambda^{\frac{-c}{2C-c}}\left(\frac{2}{C}-\frac{4}{c}\right)-\frac{2}{C}\Lambda+\frac{4}{c}-1\right]\left(m_{1,0}+m_{2,0}\right).
 \end{align*}
 The last term in the previous inequality becomes to be non positive under the conditions given in \eqref{eq:critical-inequality}. This fact implies that, there is a $T_{max}<t_{min}$ such that $\beta(T_{max})$ vanishes. This is a contradiction with the fact that solutions are smooth with positive mass. It is clear also that mass concentrates at the origin.
 \fin
 
 Finally, we consider the case when one of the species does not have a degradation term.
 
 \begin{theorem}
 	\label{theorem:blow-up-case2}
 	Under the same conditions as in Lemma \ref{lemma:estimate-second-moment} assume that $c_1=c>0$ and $c_2=0$. If one of the following conditions
 	\begin{enumerate}
 		\item[a)] $\displaystyle \alpha_2^2m_{2,0}>8\pi$,
 		\item[b)] $\Lambda>1$, $\displaystyle \alpha_2^2m_{2,0}\leq8\pi$, and $\alpha_1>0$ is large enough
 	\end{enumerate}
 	hold true, then solutions present blowing-up in a finite time.
 \end{theorem}
\noindent\emph{Proof}
Again, Lemma \ref{lemma:estimate-second-moment} reads as
\begin{align}
	e^{-ct}M_1'+M_2'\leq 4\left(m_{1,0}e^{-ct}+m_{2,0}\right)-\frac{1}{2\pi}\left(\alpha_1m_{1,0}e^{-ct}+\alpha_2m_{2,0}\right)^2=:\psi(t).\label{eq:theorem-blow-up-case2-1}
\end{align}
In this case we use the notation
$$\beta(t)=\int_{\bar{D}}|x|^2\left(\rho_1(x,\tau)+\rho_2(x,\tau)\right)\ dx+c\int_{0}^{t}\int_{\bar{D}}|x|^2\rho_1(x,\tau)\ dx\ d\tau.$$
Integrating \eqref{eq:theorem-blow-up-case2-1}, we get the inequality
\begin{align}
	\beta(t)\leq&\ (\alpha_1\alpha_2m_{2,0}-4\pi)\frac{m_{1,0}e^{-ct}}{c\pi}+\frac{\alpha_1^2m_{1,0}^2}{4\pi c}e^{-2ct}+(8\pi m_{2,0}-\alpha^2m_{2,0}^2)\frac{t}{2\pi}\nonumber\\
	\ &-(\alpha_1\alpha_2m_{2,0}-4\pi)\frac{m_{1,0}}{c\pi}-\frac{\alpha_1^2m_{1,0}^2}{4\pi c}+m_{1,0}+m_{2,0}.\label{eq:theorem-blow-up-case2-2}
\end{align}
Looking at the right hand side of \eqref{eq:theorem-blow-up-case2-2}, it is clear that $\beta(t)$ vanishes in finite time under the assumption in {\it a)}. That clearly implies blow-up of solutions in finite time as we argued in the proof of Theorem \ref{theorem:blow-up-solutions-1}.

In the other case we going to minimize the right hand side of \eqref{eq:theorem-blow-up-case2-2} whose derivative is $\psi(t)$ in \eqref{eq:theorem-blow-up-case2-1}. Note that conditions $\Lambda>1$ and $\displaystyle \alpha_2^2m_{2,0}\leq8\pi$ implies that $\psi(t)$ has a root $t_{min}$ that is the minimum procured. A direct computation shows that this root satisfies the relation,
\begin{equation}
	z:=e^{-ct_{min}}=\frac{4\pi}{\alpha_1^2m_{1,0}}\left[1+\sqrt{\left(\frac{\alpha_1m_{2,0}(\alpha_1-\alpha_2)}{2\pi}+1\right)}-\frac{\alpha_1\alpha_2m_{2,0}}{4\pi}\right].\label{eq:theorem-blow-up-case2-3}
\end{equation}
Moreover, we can rewrite \eqref{eq:theorem-blow-up-case2-2} at $t=t_{min}$ in the following way
\begin{align}
\beta(t_{min})\leq&\ (\alpha_1\alpha_2m_{2,0}-4\pi)\frac{m_{1,0}}{c\pi}\left(\frac{z}{2}-1\right)+\frac{(\alpha^2m_{2,0}^2-8\pi m_{2,0})}{2\pi c}\left(\log(z)-\frac{1}{2}\right)\nonumber\\
\ &-\frac{\alpha_1^2m_{1,0}^2}{4\pi c}+m_{1,0}+m_{2,0}.\label{eq:theorem-blow-up-case2-4}
\end{align}
Finally, looking at relations \eqref{eq:theorem-blow-up-case2-3} and \eqref{eq:theorem-blow-up-case2-4} it is clear that $\beta(t)$ vanishes in a finite time for $\alpha_1$ large enough. This implies again the existence of a finite time where blow-up of solutions occurs.
\fin
We conclude this work by noting some direct consequences of the results in this sections.
\begin{remark}
	Theorem \ref{theorem:blow-up-case2} apply to the system \eqref{eq:system-origin-two-species} even if $c_1$ and $c_2$ are positive. This is because we can change any equality in \eqref{eq:system-origin-two-species} by an inequality by removing one term of type $-c_i\rho_i$. Thus any of the following conditions implies on the blowing-up of solutions in finite time of system \eqref{eq:system-origin-two-species} for $c_1,c_2>0$:
	\begin{enumerate}
		\item  $\Lambda>1$ and $\frac{2c}{C}-4\leq\Lambda^{\frac{c}{2C-c}}\left(\frac{2c}{C}\Lambda-4-c\right)$ where $c=\min\{c_1,c_2\}$ and $C=\max\{c_1,c_2\}$,
		\item $\displaystyle \alpha_1^2m_{1,0}>8\pi$,
		\item $\displaystyle \alpha_2^2m_{2,0}>8\pi$,
		\item $\Lambda>1$, $\displaystyle \alpha_1^2m_{1,0}\leq8\pi$, and $\alpha_2>0$ is large enough,
		\item $\Lambda>1$, $\displaystyle \alpha_2^2m_{2,0}\leq8\pi$, and $\alpha_1>0$ is large enough.
	\end{enumerate}
\end{remark}
\noindent\textbf{Acknowledgement:} The author acknowledges the support of ``Concurso de Proyectos de Investigación y Fondos Semilla 2018'' from the ``Universidad Católica San Pablo'', Arequipa - Perú. Project UCSP-2018-FS-05

		\textit{E-mail address}: julioguevara08@gmail.com, jcvalencia@ucsp.edu.pe
	\end{document}